\documentclass[11pt]{amsart}
\usepackage{amsmath,amsthm,amsfonts}
\usepackage{graphicx}
\newtheorem{theorem}{Theorem}[section]

\theoremstyle{definition}
\newtheorem{algorithm}[theorem]{Algorithm}

\newtheorem{proposition}[theorem]{Proposition}

\newcommand{\gp}[1]{{\left\langle #1 \right\rangle}}
\newcommand{\rb}[1]{{\left( #1 \right)}}

\newcommand{\s}{\sigma}
\newcommand{\Bi}{B_\infty}
\newcommand{\MN}{\mathbb{N}}

\title{Cryptanalysis of Shifted Conjugacy Authentication Protocol}

\author[]{Jonathan Longrigg}
\address{Department of Mathematics, University of Newcastle, Newcastle,
NE1 7RU, UK.} \email{Jonathan.Longrigg@newcastle.ac.uk}

\author[]{Alexander Ushakov}
\address{Department of Mathematics, Stevens Institute of Technology,
Hoboken, NJ 07030, USA.} \email{sasha.ushakov@gmail.com}

\begin{document}

\maketitle

\begin{abstract}
In this paper we present the first practical attack on the shifted
conjugacy-based authentication protocol proposed by P. Dehornoy in
\cite{d}. We discuss the weaknesses of that primitive and propose
ways to improve the protocol.
\end{abstract}

\section{Introduction}

Let $B_{n}$ be the group of braids on $n$ strands given by its
standard Artin presentation
\begin{displaymath}
B_n =
 \left\langle
\begin{array}{lcl}\s_1,\ldots,\s_{n-1} & \bigg{|} &
\begin{array}{ll}
 \s_i\s_j\s_i=\s_j\s_i\s_j & \textrm{if }|i-j|=1 \\
 \s_i\s_j=\s_j\s_i & \textrm{if }|i-j|>1
\end{array}
\end{array}
\right\rangle,
\end{displaymath}
and let $B_\infty$ be the group of braids on infinitely many strands
generated by an infinite family $\{\s_1, \s_2, \ldots \}$ subject
to the same relations. There are several normal forms available
for elements of $B_n$, e.g., Garside normal form \cite[Chapter 9]{Epstein}, or the Birman-Ko-Lee
normal form \cite{BKL}. For the purposes of this paper it is convenient to
define the length of an element $x \in B_n$ to be the length of its
Garside normal form and denote it by $|x|_{\Delta_n}$.

For a braid word $w = \s_{i_1}^{\varepsilon_1} \ldots
\s_{i_k}^{\varepsilon_k}$ over the group alphabet of $B_\infty$
define a braid word
$$d(w) = \s_{i_1+1}^{\varepsilon_1} \ldots \s_{i_k+1}^{\varepsilon_k}.$$
The mapping $w \mapsto d(w)$ induces a monomorphism of $B_\infty$
and is referred to as a {\em shift operator}. Now, for braids $a,b
\in B_\infty$ define a braid $a*b$ by
\begin{equation}\label{eq:shifted_conj}
a*b = a \cdot d(b) \cdot \s_1 \cdot d(a^{-1}).
\end{equation}
The operator $*:B_\infty \times B_\infty \rightarrow
B_\infty$ is called the {\em shifted conjugacy} operator.

The Dehornoy authentication protocol is the following sequence of
steps. First, Alice prepares her public and private keys. She
randomly chooses elements $s,p \in B_n$,
and computes $p' = s \ast p$. The element $s$ is called her {\em
private key} (to be kept secret) and the pair $(p,p')$ is called
her {\em public key} (to be published).

The protocol is a Fiat-Shamir-style \cite{FS} authentication protocol
in which a single round of the protocol is performed as follows:
\begin{enumerate}
    \item[A.]
Alice chooses a random $r\in B_n$ and sends a pair $(x,x')$
(called the {\em commitment}) to Bob, where $x = r*p$ and $x' = r*p'$.
    \item[B.]
Bob chooses a random bit $b$ (called the {\em challenge}) and sends it
to Alice:
\begin{enumerate}
    \item[0)]
If $b=0$ then Alice sends $y=r$ to Bob, and Bob checks that the
equalities $x=y*p$ and $x'=y*p'$ are satisfied.
    \item[1)]
If $b=1$ then Alice sends $y=r*s$ to Bob, and Bob checks that the
equality $x'=y*x$ is satisfied.
\end{enumerate}
\end{enumerate}
Noting that $u*(v*w)=(u*v)*(u*w)$, it is straightforward to check that a correct answer $y$ to a
challenge leads to a successful check.

To break the system it is sufficient to find any $s' \in B_n$
satisfying $s' \ast p = p'$. Hence the security of this protocol
is, in particular, based on the difficulty of the Shifted
Conjugacy Search Problem (ShCSP) which is the following
algorithmic question:
\begin{quote}
For a pair of braids $p,p' \in \Bi$ find a braid $s \in \Bi$ such
that $p' = s * p$ (provided that such $s$ exists).
\end{quote}
Similarly one can formulate the Shifted Conjugacy Decision Problem (ShCDP):
\begin{quote}
For a pair of braids $p,p' \in \Bi$ determine if there exists a
braid $s \in \Bi$ such that $p' = s * p$.
\end{quote}
These problems first appeared in \cite{d} and were not explored
enough to give a precise answer about their time complexity or
even decidability (for the decision problem). Despite the
resemblance to the conjugacy decision problem (CDP) which is
decidable and suspected to have polynomial-time solution (see
\cite{BGM1}, \cite{BGM2}), it is not clear if ShCDP is solvable.

It is not discussed in Dehornoy's original paper \cite{d} how to generate public
and private keys, the author just proposes a primitive and
provides some intuition on why the primitive might be hard. In
this paper we use the easiest method of key generation: keys are chosen
uniformly from the ambient free group and then considered as
words in the braid group.
\begin{enumerate}
    \item[1)]
Fix the rank of the braid group $n$. It is an important parameter
in the scheme, efficiency of all the operations depends on it.
    \item[2)]
Fix numbers $L$ and $K$, the key lengths. These numbers are the main security
parameters. (In all our experiments $L=K$.)
    \item[3)]
Pick randomly and uniformly a braid word $p$ (resp. $s$) from the
set of all braid words of length $L$ (resp. $K$.)
    \item[4)]
Finally, compute $p' = s \ast p$.
\end{enumerate}

To summarize the results of our work:
\begin{enumerate}
    \item[A.]
Even though it seems unlikely that ShCSP can be deterministically
reduced to CSP we argue that ShCSP can be reduced to CSP {\em
generically} (for most of the inputs) and present the reduction.
    \item[B.]
We present the results of actual experiments. The following table shows the
percentage of success in our experiment. The number of
successes where the result was equal to the original key is in brackets.
For instance, if the key length is $100$ and the platform group is $B_{40}$
then, using our attack, in 24\% of cases we recovered an element
$s\in B_n$ such that $p'=s*p$ and hence broke the protocol.
In $10\%$ of the cases the obtained element was equal to the original $s\in B_n$
generated by Alice.
\begin{table}[h]
\begin{center}
\begin{tabular}{|l|r|r|r|}
\hline
Key length, $K=L$           & 100 & 400 & 800 \\
\hline \hline
$B_{10}$                  & 100(100) & 100(100) & 100(100) \\
\hline
$B_{40}$                  & 24(10)  & 99(99)  & 92(92)  \\
\hline
$B_{80}$                  & 2(0)   & 47(39)  & 92(92) \\
\hline
\end{tabular}
\end{center}
\caption{\label{tb:success} Success rate in our experiments.}
\end{table}
    \item[C.]
We analyze the results and make several recommendations on how to
generate hard instances of ShCSP.
\end{enumerate}

The paper is organized as follows. In Section \ref{se:attack} we
describe our heuristic algorithm and argue that it works for most
inputs. In Section \ref{se:exper_results} we present more
detailed results of experiments (than in the table above), discuss
the reasons for success/failure, and make suggestions on the
generation of hard keys.

All the algorithms described in this paper are available at
\cite{CRAG}.



\section{The attack}
\label{se:attack}

In this section we present the mathematical background for the attack.
For $n \in \MN$ define braids
    $$\delta_{n} = \sigma_{n-1} \ldots \sigma_1,$$
    $$\Delta_{n} = (\sigma_{n-1} \ldots \sigma_1) \cdot (\sigma_{n-1} \ldots \sigma_2) \cdot \ldots  \cdot (\sigma_{n-1}).$$
It is easy to check that for any $i = 1,\ldots,n-1$ the following
equality holds in $B_{n+1}$:
\begin{equation}\label{eq:equat0}
\delta_{n+1}^{-1} \sigma_i \delta_{n+1} =_{B_{n+1}} \sigma_{i+1} =
d(\s_{i}).
\end{equation}

\begin{proposition} \label{pr:problem_equivalence}
Let $p, p', s \in B_{n}$. Then $s$ satisfies the shifted conjugacy
equation for $p$ and $p'$
\begin{equation}\label{eq:equat1}
    p' =_{B_{n+1}} s*p,
\end{equation}
if and only if it satisfies the conjugacy equation for $p'
\delta_{n+1}^{-1}$ and $d(p) \sigma_1 \delta_{n+1}^{-1}$
\begin{equation}\label{eq:equat2}
    p' \delta_{n+1}^{-1} =_{B_{n+1}} s  \cdot d(p) \sigma_1 \delta_{n+1}^{-1} \cdot s^{-1}.
\end{equation}
\end{proposition}

\begin{proof}
Follows from (\ref{eq:equat0}).

\end{proof}

\begin{proposition} \label{pr:problem_solving}
Let $p, p', s \in B_{n}$ be braids satisfying $p'=_{B_{n+1}}s*p$,
and let $s'\in B_{n+1}$. Then
    $$p' \delta_{n+1}^{-1}=_{B_{n+1}} s'  \cdot d(p) \sigma_1 \delta_{n+1}^{-1} \cdot s'^{-1} \iff s'^{-1} s\in C_{B_{n+1}}\rb{d(p) \s_1 \delta_{n+1}^{-1}}$$
where $C_{B_{n+1}}\rb{d(p) \sigma_1 \delta_{n+1}^{-1}}$ denotes
the centralizer of $d(p) \sigma_1 \delta_{n+1}^{-1}$ in $B_{n+1}$.
\end{proposition}

\begin{proof}
Obvious.

\end{proof}

Now, it follows from Propositions \ref{pr:problem_equivalence} and
\ref{pr:problem_solving} that the ShCSP can be solved in two
steps:
\begin{enumerate}
    \item[(S1)]
Find a solution $s' \in B_{n+1}$ of the equation $p'
\delta_{n+1}^{-1}=_{B_{n+1}} s' \cdot d(p) \sigma_1
\delta_{n+1}^{-1} \cdot s'^{-1}$. This can be done using the ultra
summit set technique invented in \cite{G}.
    \item[(S2)]
``Correct'' the element $s' \in B_{n+1}$ to obtain a solution
$s\in B_n$ of (\ref{eq:equat1}), i.e., find a suitable element $c
\in C_{B_{n+1}} \rb{d(p) \sigma_1 \delta_{n+1}^{-1}}$ such that $t
= s' c\in B_n$ satisfies (\ref{eq:equat1}). We refer to this step
as a {\em centralizer attack}.
\end{enumerate}
The description of ultra summit sets is out of the scope of our
paper, so we omit details on step (S1). Step (S2) requires some
elaboration. To be able to work with elements of $C = C_{B_{n+1}}
\rb{d(p) \sigma_1 \delta_{n+1}^{-1}}$ efficiently we need to
describe $C$ in some convenient way, for instance, by providing a
set of generators. Hence, step (S2) itself consists of two smaller
steps: capturing $C$ and finding the required element $c \in C$.

The only known algorithm \cite{GM-centr} for computing a
generating set for a centralizer reduces to the construction of
so-called super summit sets, the size of which is not known to be
polynomially-bounded, and which is usually hard in practice.
Hence, the approach of describing the whole generating set does
not seem feasible. Instead, we can work with the subgroup of
$C_{B_{n+1}} \rb{d(p) \sigma_1 \delta_{n+1}^{-1}}$ considered in
Proposition \ref{pr:subgp_C}. For $p \in B_n$ define braids:
    $$c_1 = \Delta_{n+1}^2,~~ c_2 = d(p) \sigma_2^{-1} \ldots \sigma_n^{-1},~~ c_3 = \sigma_1 \ldots \sigma_n^2 \sigma_{n-1}\ldots \sigma_{1},$$
and
    $$c_4 = c_1^{-1},~ c_5 = c_2^{-1},~ c_6 = c_3^{-1}.$$

\begin{proposition}\label{pr:subgp_C}
Let $p \in B_n$ and $C = C_{B_{n+1}} \rb{d(p) \sigma_1
\delta_{n+1}^{-1}}$. The following holds:
\begin{itemize}
    \item
$c_1,c_2, c_3 \in C$,
    \item
$C' = \gp{c_1,c_2,c_3}$ is an abelian subgroup of $B_{n+1}$ and,
hence, has polynomial growth.
\end{itemize}
\end{proposition}

\begin{proof}
Observe that the equality
    $$d(p) \sigma_1 \delta_{n+1}^{-1} =_{B_{n+1}} d(p) \sigma_2^{-1} \ldots \sigma_n^{-1}$$
holds in $B_{n+1}$ and the element $d(p) \sigma_2^{-1} \ldots
\sigma_n^{-1}$ involves generators $\sigma_2, \ldots ,\sigma_n$
only. It is intuitively obvious that $c_2$ commutes with $c_3$
when you observe that in the braid diagram for $c_3$ none of
strands $2$ to $n$ cross over. Furthermore, $c_1$ generates the
center of $B_{n+1}$. Thus, the subgroup $C'$ is abelian.

\end{proof}

Now, having fixed the subgroup $C' = \gp{c_1, c_2, c_3}$ we can
describe the heuristic procedure for finding the required $c \in
C'$. For any braid $t \in B_n$ define
    $$l_t = |p'^{-1} (t \ast p) |_{\Delta_{n+1}}$$
 and observe that $t = s' c$ satisfies $t \ast p =
p'$ if and only if $l_t =  0$ and the value of $l_t$ can be used
to guide our heuristic search, the smaller $l_t$ the ``closer" $t$
to the actual solution.

We summarize the ideas of this section into an heuristic algorithm
(Algorithm \ref{al:attack}) which for a pair of braids $p,p' \in
B_n$ attempts to find $s \in B_n$ such that $p' =_{B_{n+1}} s\ast
p$. The algorithm starts out by finding a solution $s'$ to the
conjugacy equation (\ref{eq:equat2}).
    It keeps two sets: $S$ (elements in working) and $M$ (worked out
elements) of pairs $(t,l_t)$ where $t \in B_{n+1}$ is a possible
solution and $l_t = |p'^{-1} (t \ast p) |_{\Delta_{n+1}}$.
    Initially, we have $S = \{(s',l_{s'})\}$, where $s'$ is found in step
(S1) and $l_{s'} = |p'^{-1} (s'\ast p) |_{\Delta_{n+1}}$.
    On each iteration we choose a pair $(t,l_t)$ from $S$ with the
smallest value $l_t$ (the ``fittest" one). If $l_t = 0$ then $t$
is a solution.
    If $l_t \ne 0$ then compute new possible solutions $t_i = t c_i$
and add corresponding pairs $(t_i,l_{t_i})$ into $S$. After all
$(t_i,l_{t_i})$ are added to $S$ the current pair $(t,l_t)$
becomes worked out.

\begin{algorithm}\label{al:attack}
({\em Heuristic algorithm for solving ShCSP})
    \\{\sc Input:}
Braids $p,p' \in B_n$.
    \\{\sc Output:}
A braid $s\in B_n$ such that $p' =_{B_\infty} s \ast p$.
    \\{\sc Computation:}
\begin{enumerate}
    \item[A.]
Using the ultra summit set technique compute $s' \in B_{n+1}$
satisfying $p' \delta_{n+1}^{-1} =_{B_{n+1}} s' \cdot d(p)
\sigma_1 \delta_{n+1}^{-1} \cdot s'^{-1}$.
    \item[B.]
Put $S = {(s',~|p'^{-1} (s' \ast p)|_{\Delta_{n+1}})}$ and $M =
\emptyset$.
    \item[C.]
Until a solution is found:
\begin{enumerate}
    \item[1.]
Choose a pair $(t,l_t)$ from $S$ with the smallest $l_t$.
    \item[2.]
If $l_t=0$ then output $t$.
    \item[3.]
Otherwise for each $i=1, \ldots, 6$
\begin{enumerate}
    \item
Compute $t_i = t \cdot c_i$ and $l_{t_i} = |p'^{-1} (t_i \ast
p)|_{\Delta_{i+1}}$
    \item
If $(t_i,l_{t_i})$ belongs neither to $S$ nor to $M$ then add it
into $S$.
\end{enumerate}
    \item[4.]
Remove the current pair $(t,l_t)$ from $S$ and add it into $M$.
\end{enumerate}
\end{enumerate}
\end{algorithm}

\section{Experimental results}
\label{se:exper_results}

Algorithm \ref{al:attack} always produces the correct answer when it
halts, though it does not always stop. There are two possible
reasons for a failure.
\begin{enumerate}
    \item[1)]
Failure on step A. The precise complexity of the ultra summit set
algorithm is not known, though it is proved that for certain
classes of braids it is polynomial. But even if it is polynomial,
the degree of a polynomial can be too large to be used in
practical computations.
    \item[2)]
Failure in the loop C. There are two possible reasons for this.
The first reason is that we use a subgroup $C'$ of $C =
C_{B_{n+1}} \rb{d(p) \sigma_1 \delta_{n+1}^{-1}}$ which can be a
proper subgroup and it might not contain the required element.

The second reason is that the heuristic based on the length $|p'^{-1} (s_i \ast
p)|_{\Delta_{n+1}}$ might be bad.
\end{enumerate}

To test the efficiency of the algorithm we performed a series of
experiments in which we limited the time allowed for each part:
$1$ hour for step A and 15 minutes for step C on a personal
computer (CPU 2.66 GHz). The percentage of success in experiments for different
parameters is shown in Table \ref{tb:success}. Below we present
detailed information on each step of computations.

Recall that in step A we solve the conjugacy search problem for
braids $p' \delta_{n+1}^{-1}$ and $d(p) \sigma_1
\delta_{n+1}^{-1}$ in $B_{n+1}$. The percentage of failure on step
A is shown in Table \ref{tb:failure_A}. We see that the shorter
the length of the key $p$ relative to the rank of the braid group
the harder it becomes to solve the conjugacy search problem. We
cannot explain this phenomenon, the method of ultra summit sets is
very difficult to analyze. We suspect that the reason for such
behavior is that shorter words are less likely to be rigid (expose
free-like behavior when cycling.)
\begin{table}[h]
\begin{center}
\begin{tabular}{|l||r|r|r|}
\hline
Key length                & 100 & 400 & 800 \\
\hline \hline
$B_{10}$                  &   $0$ of $100$ ($0.00\%$)  &   $0$ of $100$ ($0.00\%$)  &  $0$ of $100$ ($0.00\%$) \\
\hline
$B_{40}$                  &   $9$ of $100$ ($9.00\%$)  &   $0$ of $100$ ($0.00\%$)  &  $1$ of $100$ ($0.00\%$)\\
\hline
$B_{80}$                  &  $59$ of $100$ ($59.00\%$) &  $12$ of $100$ ($12.00\%$) &  $6$  of $100$ ($6.00\%$) \\
\hline
\end{tabular}
\end{center}
\caption{\label{tb:failure_A} Failure on step A}
\end{table}

The percentage of failure in loop C of Algorithm \ref{al:attack},
given that step A was successfully completed, is shown in Table
\ref{tb:failure_C}.  Essentially, we see the same pattern of failure
as in the solution of the conjugacy problem in Table
\ref{tb:failure_A}. The shorter the length of the key relative to
the rank of the braid group the harder it is to find a suitable
element of the centralizer. These results are easier to explain.
Basically the shorter the key the bigger the centralizer. Hence,
for shorter elements it is more likely that $C'$ is a proper
subgroup of $C = C_{B_{n+1}} \rb{d(p) \sigma_1 \delta_{n+1}^{-1}}$
and less likely that $C'$ actually contains the required element
$c$.
\begin{table}[h]
\begin{center}
\begin{tabular}{|l||r|r|r|}
\hline
Key length                & 100 & 400 & 800 \\
\hline \hline
$B_{10}$                  &   $0$ of $100$ ($0.00\%$) &   $0$ of $100$ ($0.00\%$) &   $0$ of $100$ ($0.00\%$) \\
\hline
$B_{40}$                  &  $67$ of $91$ ($73.62\%$) &   $1$ of $100$ ($1.00\%$) &   $7$ of $99$ ($7.07\%$) \\
\hline
$B_{80}$                  &  $39$ of $41$ ($95.12\%$) &  $41$ of $88$ ($46.59\%$) &   $2$ of $94$ ($2.12\%$) \\
\hline
\end{tabular}
\end{center}
\caption{\label{tb:failure_C} Failure in loop C}
\end{table}

\section{Conclusions}

In this section we discuss methods of key generation
invulnerable to the attack proposed in Section \ref{se:attack}.
Recall that the success of the attack relies on two properties of
braids $p$, $p'$:
\begin{enumerate}
    \item
the conjugacy search problem is easy for the pair $(p'
\delta_{n+1}^{-1}, d(p) \sigma_1 \delta_{n+1}^{-1})$ in $B_{n+1}$;
    \item
the centralizer $C_{B_{n+1}}(d(p) \sigma_1 \delta_{n+1}^{-1})$ is
``small" (isomorphic to an Abelian group of small rank.)
\end{enumerate}
If either of the two properties above is not satisfied then the attack
is likely to fail. Even though there is no known polynomial time algorithm
solving CSP for braids, recent developments of \cite{BGM1} and
\cite{BGM2} suggest that CSP might be easy for generic braids
(pseudo-Anosov type) and it might be difficult to randomly
construct hard instances for CSP. Though there is no proof yet
that the pseudo-Anosov type of braids is generic it is a common belief
that this is so. Another interesting recent development is
\cite{KLT} where the authors present a few braids with very
large ultra summit sets.

The only part which can be controlled is the growth and structure
of the centralizer $C$. As we mentioned above the only known
algorithm for computing a generating set for a centralizer reduces
to construction of super summit sets which is not known to be
polynomially hard and is usually very inefficient practically. We
can choose $p$ so that $C(d(p) \sigma_1 \delta_{n+1}^{-1})$ is a
large non-Abelian group. For more on the structure of centralizers
see \cite{GMW}. For ideas on how to generate elements with large
centralizers see \cite{su}.

\end{document}